# A polynomial time knot polynomial

Dror Bar-Natan and Roland van der Veen

**Abstract**

We present the strongest known knot invariant that can be computed effectively (in polynomial time)[1].

## 1 Introduction

We present a simple, strong knot invariant that is closely related to the Alexander polynomial and seems to share many of its good properties. For example, unlike the commonly used quantum invariants such as the Jones polynomial, our invariant is computable in polynomial time. For simple knots the Jones polynomial works well but as the number of crossings grows the exponentially many terms in the resulting state sums quickly become unmanageable. Our invariant does not suffer from these issues, it scales well with the complexity of the knot.

The plan of the paper is as follows. Our invariant is based on normal ordering exponentials in the $q$-Weyl algebra $\langle E, F \rangle / (EF - qFE = 1)$, see Section 5. As a warm-up we first show how the Alexander polynomial may be derived from the ordinary Weyl algebra $\langle E, F \rangle / (EF - FE = 1)$ in Section 4. These algebras are but examples of a more general theory of invariants coming from algebras satisfying certain equations that we introduce in Section 2. As a preview we start by giving a condensed definition of the knot invariant. Proofs and additional (conjectured) properties are discussed in Section 6.

### 1.1 The knot invariant

Consider a (long) knot $K$ presented as a proper smooth embedding of $[0, 1]$ into the closed unit ball such that the projection on the third coordinate is a generic immersion $\gamma$ in the plane, see for example Figure 1. More specifically, assume that there is an $n \in \mathbb{N}$ such that $\gamma$ has the following properties. The points $\gamma(\frac{k}{n+1})$ where $k \in \{1, \ldots, n\}$ are the union of all double points and all points where $\gamma'$ is in the direction of the positive $x$-axis. The double points are known as crossings and the latter as cuaps. Close to any crossing we assume $\gamma'$ has positive $y$-coordinate. The sign of a crossing is the sign of the $x$-coordinate of $\gamma'$ at the overpass. A crossing is denoted $X^\sigma_{i,j}$ where $\sigma$ is the sign and $i, j$ are the labels of the over and under strand. The sign of a cuap is the sign of the $y$-direction of $\gamma''$. A cuap is denoted $u^\sigma_i$ where $\sigma$ is the sign and $i$ is its label.

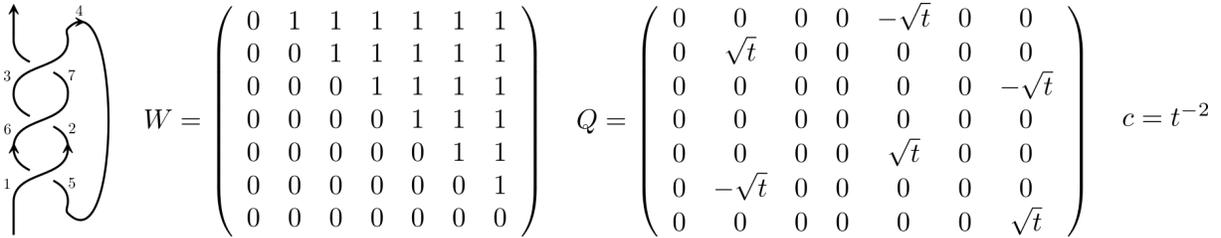

Figure 1: A diagram for the trefoil knot $3_1$. The double points at the crossings and the right-pointing cuaps are enumerated in order of appearance. The matrices $W, Q$ and the number $c$ necessary for computation of the invariant $Z_0$ are listed next to it.

To define our knot invariant $Z_0$ we need to introduce some preliminary constructions. Let $E^i_j$ be the elementary $n \times n$ matrix with a single non-zero entry 1 at the $(i, j)$-th place. Define the matrices

$$Q = \sum_{X^\sigma_{i,j}} \sigma t^{\frac{\sigma}{2}}(E^j_j - E^i_j) \quad W = \sum_{i<j} E^i_j \quad c = \prod_{X^\sigma_{i,j}, u^\sigma_i} t^{-\sigma}$$

Here the sum and product are over all crossings/cuaps in the diagram. Next recall the adjugate of a matrix $M$ is defined by $M \operatorname{adj}(M) = \det(M)I$ and define

$$B = I - (t^{\frac{1}{2}} - t^{-\frac{1}{2}})WQ \quad G = Q\operatorname{adj}(B) \quad H = \operatorname{adj}(B)W$$

---

[1] This work was partially supported by NSERC grant RGPIN 262178 and the Netherlands Organisation for Scientific Research.

$$Z_G = (t - t^{-1}) \sum_{j=2}^{n} \sum_{a,b<j} G_a^j \left( \frac{1}{2} G_b^j + \sum_{g>j} G_b^g \right)$$

$$Z_H = \sum_{X_{i,j}^\sigma} \frac{\sigma}{2}((1-t^\sigma)H_i^j)^2 - \frac{\sigma}{2}((1+t^\sigma)H_j^j)^2 + \sigma t^\sigma (H_i^j H_j^i + H_i^i H_j^j) + t^\sigma (1-t) H_i^j ((1+\sigma)H_j^j + (1-\sigma)H_i^i)$$

$$+ \det(B) \sum_{u_i^\sigma} \sigma H_i^i$$

**Theorem 1.** $c^{\frac{1}{2}} \det(B)$ *is the Alexander polynomial* $\Delta$ *and* $Z_0 = c(Z_G + Z_H)$ *is a knot invariant. Both are elements of* $\mathbb{Z}[t, t^{-1}]$ *computable in* **polynomial time**.

The above formulas immediately show that the computation must be polynomial time, for a more detailed discussion see Theorem 5 in Section 6.

The pair $\Delta, Z_0$ distinguishes all knots in the Rolfsen table of prime knots up to ten crossings, see the Appendix. That is a better performance than for example Khovanov homology and HOMFLY polynomial combined[2]. More importantly $Z_0$ appears to share many of the desirable properties of the Alexander polynomial.

Our invariant $Z_0$ appears to coincide with a part of the colored Jones invariant studied by Overbay [12] and Rozansky [14], see Conjecture 3 of Section 6. The present approach seems simpler as it allows a local version and a polynomial time algorithm. We also conjecture new bounds on the knot genus and argue $Z_0$ detects mirror images.

## 1.2 Example: Trefoil

We illustrate the computation of $Z_0$ for the trefoil knot $3_1$ shown in Figure 1.

$$B = \begin{pmatrix} 1 & 0 & 0 & 0 & 1-t & 0 & 0 \\ 0 & t & 0 & 0 & 1-t & 0 & 0 \\ 0 & t-1 & 1 & 0 & 1-t & 0 & 1-t \\ 0 & t-1 & 0 & 1 & 1-t & 0 & 1-t \\ 0 & t-1 & 0 & 0 & 1 & 0 & 1-t \\ 0 & 0 & 0 & 0 & 0 & 1 & 1-t \\ 0 & 0 & 0 & 0 & 0 & 0 & 1 \end{pmatrix} \quad \Delta_{3_1}(t) = c^{\frac{1}{2}} \det(A) = t - 1 + t^{-1}$$

$$G = \begin{pmatrix} 0 & t^{\frac{3}{2}} - t^{\frac{1}{2}} & 0 & 0 & -t^{\frac{3}{2}} & 0 & t^{\frac{3}{2}} - t^{\frac{5}{2}} \\ 0 & t^{\frac{1}{2}} & 0 & 0 & t^{\frac{3}{2}} - t^{\frac{1}{2}} & 0 & t^{\frac{5}{2}} - 2t^{\frac{3}{2}} + t^{\frac{1}{2}} \\ 0 & 0 & 0 & 0 & 0 & 0 & -t^{\frac{5}{2}} + t^{\frac{3}{2}} - t^{\frac{1}{2}} \\ 0 & 0 & 0 & 0 & 0 & 0 & 0 \\ 0 & t^{\frac{1}{2}} - t^{\frac{3}{2}} & 0 & 0 & t^{\frac{3}{2}} & 0 & t^{\frac{5}{2}} - t^{\frac{3}{2}} \\ 0 & -t^{\frac{1}{2}} & 0 & 0 & t^{\frac{1}{2}} - t^{\frac{3}{2}} & 0 & -t^{\frac{5}{2}} + 2t^{\frac{3}{2}} - t^{\frac{1}{2}} \\ 0 & 0 & 0 & 0 & 0 & 0 & t^{\frac{5}{2}} - t^{\frac{3}{2}} + t^{\frac{1}{2}} \end{pmatrix} \quad Z_G = t^4 - \frac{3t^2}{2} + \frac{1}{2}$$

$$H = \begin{pmatrix} 0 & t^2 - t + 1 & t & t & t & t^2 & t^2 \\ 0 & 0 & 1 & 1 & 1 & t & t \\ 0 & 0 & t - t^2 & 1 & 1 & t & t \\ 0 & 0 & t - t^2 & t - t^2 & 1 & t & t \\ 0 & 0 & 1 - t & 1 - t & 1 - t & 1 & 1 \\ 0 & 0 & 0 & 0 & 0 & 0 & t^2 - t + 1 \\ 0 & 0 & 0 & 0 & 0 & 0 & 0 \end{pmatrix} \quad Z_H = t^4 - 3t^3 + \frac{7t^2}{2} - t - \frac{1}{2}$$

It follows that $Z_0(3_1) = c(Z_G + Z_H) = 2 - t^{-1} - 3t + 2t^2$ and its normalization is $\rho_1(t) = -t - t^{-1}$, see section 6. Comparing to the table in the Appendix, notice we used the mirror image of the usual trefoil $3_1$ from the table. Conjecture 1 explains the ensuing sign in $\rho_1$.

## 2 Snarl diagrams: A local version of the knot invariant

So far we defined $Z_0$ for knot diagrams but did not yet show its value is independent of the chosen diagram. As any two diagrams for the same knot are related by Reidemeister moves, all we need to do is show $Z_0$ is unchanged under those moves. Instead of attempting a direct proof we first extend $Z_0$ to a function $Z$ on more general diagrams that we call snarl diagrams. Showing invariance of $Z$ is vastly simpler since now all computations become local. Our treatment is closely related to Kauffman's rotational virtual tangles [7] but avoids virtual crossings.

---

[2] The knots $8_{16}$ and $10_{156}$ have identical Khovanov homology and identical HOMFLY polynomials [1].

**Definition 1.** *A* **snarl**[3] *diagram is a finite set $L$ together with a finite oriented graph $G = (V, E)$ and functions $\sigma : V \to \{\pm 1\}$ and $\rho : E \to \mathbb{Z}$. The edges $E$ are assumed to be a disjoint union of oriented paths and each path is labelled by an element of $L$. Furthermore the edges around any vertex are ordered cyclically such that two adjacent edges enter and two exit each vertex that is not an endpoint of a path.*

The vertices of the graph should be viewed as the crossings and endpoints of a projection of a piece of a knot. The paths labeled by $L$ correspond to the connected components. The map $\rho$ keeps track of rotation numbers of the tangent vector on the edges so that our diagrams are like Morse diagrams. To build a snarl diagram from any knot diagram just make the tangent vector near each crossing point upwards and count the resulting rotation numbers of the tangent vector at each edge. For example the edge labeled 4 in Figure 1 has $\rho = -1$ as it rotates clockwise. $\rho = 0$ for all other edges. A more interesting example of a snarl diagram is shown in Figure 2.

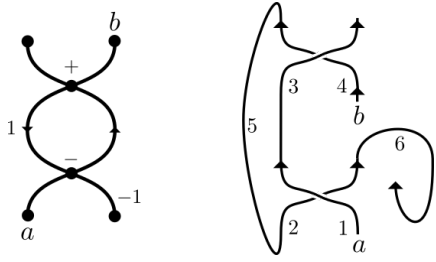

Figure 2: Left: A snarl diagram corresponding to the left-hand side of equation (2), the numbers give the value of $\rho$, the letters $a, b$ are the labels of the two components. Right: A piece of a knot whose snarl diagram would be on the left hand side. Again $a, b$ depict the labels of the two components and the integers are labels of the smaller pieces used to build up this diagram as in equation (2).

Knot diagrams may be assembled from simpler pieces by the following two operations on snarl diagrams.

**Definition 2. Disjoint union:** *For two snarl-diagrams $G, G'$ with label sets $L, L'$ the disjoint union $G \sqcup G'$ is the snarl diagram with underlying graph as indicated and label set $L \sqcup L'$. To avoid clutter we sometimes omit the $\sqcup$ symbol and use juxtaposition instead.*

**Stitching:** *For $i \neq j \in L$ and $k \notin L - \{i, j\}$ define the snarl diagram $m_k^{ij}(G)$ to be the graph obtained from $G$ by connecting the endpoint of component $i$ to the start of component $j$, erasing the vertex in the middle. For the newly created edge $e$ we define $\rho(e)$ to be the sum of the values of $\rho$ on the edges that disappear. The newly created component is labeled $k$ so the label set is $L - \{i, j\} \cup \{k\}$.*

Notice how any snarl diagram may be constructed using disjoint union and stitching from two types of fundamental diagrams: The crossing $X_{ij}^{\pm}$ where we label the over-strand $i$ and $\rho$ of every edge is 0 and the diagram $\alpha_i^r$, a single oriented edge labelled $i$ with rotation number $\rho = r$. It is sometimes useful to stitch many ends at the same time. For a sequence $I = (I_1, \ldots, I_n)$ of $n$ distinct elements of $L$ and $k \notin L - I$ define[4]
$$m_k^I = m_k^{I_1, I_2} /\!\!/ m_k^{k, I_3} /\!\!/ \ldots m_k^{k, I_{n-1}} /\!\!/ m_k^{k, I_n}$$
Even more generally if $\tau = (\tau^1, \ldots \tau^b)$ is a sequence of $b$ such sequences whose disjoint union is $L$ and $B = (B_1, \ldots, B_b)$ is any $b$-element sequence then define $m_B^\tau = m_{B_1}^{\tau^1} /\!\!/ \ldots /\!\!/ m_{B_b}^{\tau^b}$. When $B$ is not specified we take it to be $(1, \ldots, b)$, the commas are sometimes dropped for brevity.

Snarl diagrams are meant to generalize Morse diagrams of pieces of knots. As such we should consider them up to equivalence under a version of the Reidemeister moves as for example described in Chapter 3 of [11]. Using rotation numbers instead of cuaps and using [13] we only have to look at the four equations below, depicted in Figure 3.

**Definition 3.** *Consider the equivalence relation $\sim$ on the set of snarl diagrams generated by relabeling components and the equivalences*

$$X_{13}^{\pm} \sqcup \alpha_2^{\mp 1} /\!\!/ m^{(123)} \sim \alpha_1^0 \sim X_{31}^{\pm} \sqcup \alpha_2^{\pm 1} /\!\!/ m^{(123)} \qquad (1)$$

$$X_{12}^{-} \sqcup X_{34}^{+} \sqcup \alpha_5 \sqcup \alpha_6^{-1} /\!\!/ m^{(13)(4526)} \sim \alpha_1^0 \sqcup \alpha_2^0 \sim X_{12}^{+} \sqcup X_{34}^{-} \sqcup \alpha_5 \sqcup \alpha_6^{-1} /\!\!/ m^{(5163)(42)} \qquad (2)$$

$$X_{12}^{\pm} \sqcup X_{34}^{\pm} \sqcup X_{56}^{\pm} /\!\!/ m^{(13)(25)(46)} \sim X_{12}^{\pm} \sqcup X_{34}^{\pm} \sqcup X_{56}^{\pm} /\!\!/ m^{(35)(16)(24)} \qquad (3)$$

$$X_{12}^{\pm} \sim \alpha_1^{\pm 1} \sqcup \alpha_2^{\pm 1} \sqcup X_{34}^{\pm} \sqcup \alpha_5^{\mp 1} \sqcup \alpha_6^{\mp 1} /\!\!/ m^{(135)(246)} \qquad (4)$$

Our motivation for defining this equivalence relation is the following Reidemeister type theorem.

**Theorem 2.** *Every isotopy class of long knots can be represented by a snarl diagram. Two snarl diagrams representing isotopic knots are equivalent.*

---
[3]Dictionary entry: a knot or tangle, also a growl.
[4]In what follows $f /\!\!/ g$ means the composition $g \circ f$.

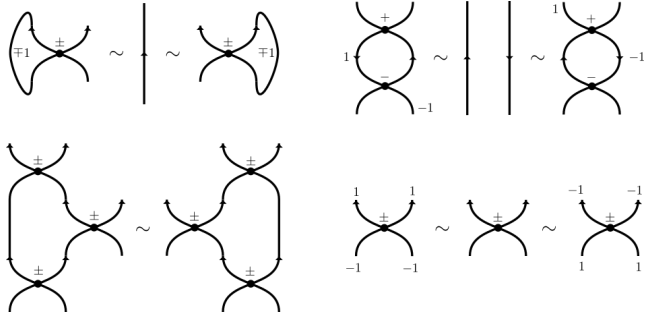

Figure 3: Equivalences on snarl diagrams. The labels of the components and degree 1 vertices are not shown. Only the edges with a non-zero value of $\rho$ are marked.

*Proof.* This can be proved exactly as in the case of the standard Reidemeister theorem for Morse diagrams of tangles. See for example Chapter 3 of [11] or [15]. The fact that our set of moves is sufficient follows from Thm 1.2 of [13]. □

Snarl diagrams may also be used to represent more general knotted objects such as planar tangles and for those a similar theorem will hold.

## 3 Snarl algebra

In this section we formulate conditions on an algebra to yield a knot invariant. Such invariants are sometimes known as universal invariants, see [6] and the references therein.

By an algebra $\mathcal{A}$ we mean a ring with 1 whose center includes a commutative base ring $R$ with 1. All tensor products are over $R$. Similar to the labeling of components of snarl diagrams we would like to explicitly label the tensor factors in our tensor products. Given a finite set $S$ we define $\mathcal{A}^{\otimes S}$ to be the tensor product of $|S|$ copies of $\mathcal{A}$ indexed by the elements of $S$. When $T \subset S$ define $\iota_T : \mathcal{A}^{\otimes T} \to \mathcal{A}^{\otimes S}$ by the identity on the tensor factors indexed by elements of $T$ and setting all other tensor factors to 1. We use the shorthand $\iota_{\{i,j\}}(Y) = Y_{ij}$ and $\iota_{\{i\}}(y) = y_i$. Denote by $m_k^{ij}$ the multiplication operation $\mathcal{A}^{\otimes S \sqcup \{i,j\}} \to \mathcal{A}^{\otimes S \sqcup \{k\}}$ defined by deleting the $i$-th and $j$-th tensor factors and placing the product of those elements into the $k$-th tensor factor.

**Definition 4.** A **snarl-algebra** is an algebra $\mathcal{A}$ together with invertible elements $X \in \mathcal{A}^{\otimes \{1,2\}}, \alpha \in \mathcal{A}$ such that the equations shown below are satisfied. For any snarl diagram $D$ with label set $L$ denote by $Z(D) \in A^{\otimes L}$ the unique element characterized by $Z(X_{ij}^\pm) = X_{ij}^\pm$ and $Z(\alpha_i^r) = \alpha_i^r$ and $Z(D \sqcup D') = Z(D) \otimes Z(D')$ and $Z(m_k^{ij} D) = m_k^{ij} Z(D)$. $Z = Z_{\mathcal{A},X,\alpha}$ is called the snarl invariant corresponding to $\mathcal{A}$.

$$Z(X_{13}^\pm \alpha_2^{\mp 1} /\!/ m^{(123)}) = Z(\alpha_1^0) = Z(X_{31}^\pm \alpha_2^{\pm 1} /\!/ m^{(123)}) \tag{5}$$

$$Z(X_{12}^- X_{34}^+ \alpha_5 \alpha_6^{-1} /\!/ m^{(13)(4526)}) = Z(\alpha_1^0 \alpha_2^0) = Z(X_{12}^+ X_{34}^- \alpha_5 \alpha_6^{-1} /\!/ m^{(5163)(42)}) \tag{6}$$

$$Z(X_{12}^\pm X_{34}^\pm X_{56}^\pm /\!/ m^{(13)(25)(46)}) = Z(X_{12}^\pm X_{34}^\pm X_{56}^\pm /\!/ m^{(35)(16)(24)}) \tag{7}$$

$$Z(X_{12}^\pm) = Z(\alpha_1^{\pm 1} \alpha_2^{\pm 1} X_{34}^\pm \alpha_5^{\mp 1} \alpha_6^{\mp 1} /\!/ m^{(135)(246)}) \tag{8}$$

This definition is designed to make following theorem hold.

**Theorem 3.** *For any snarl algebra the corresponding invariant $Z$ is well-defined and its value is independent of the snarl diagram chosen to compute it.*

*Proof.* For the well-definedness of $Z$ we argue that its value does not depend on the order of stitchings and disjoint unions used to build up the snarl. First one may carry out all disjoint union operations at the beginning. By associativity of both operations it is then clear that the order is irrelevant.

By construction equivalent snarl diagrams will yield the same value of $Z$. Theorem 2 finishes the proof. □

As a relatively simple example, not used below, consider the group algebra $\mathbb{C}G$ of a finite group $G$ and its dual $\mathbb{C}(G) = \{f : G \to \mathbb{C}\}$. Define a snarl algebra $D(G) = \mathbb{C}(G) \otimes \mathbb{C}G$ with multiplication determined by

$$(\delta^g \otimes h)(\delta^{g'} \otimes h') = \delta^g \delta^{hg'h^{-1}} \otimes hh'$$

here $\delta^g$ denotes the function that takes value 1 on $g$ and is zero otherwise.

The reader is invited to check that $D(G)$ becomes a snarl algebra if we set:

$$Z_{D(G)}(X_{ij}^\sigma) = \sum_{g \in G}(\delta^g \otimes g^{-\sigma})_i (1 \otimes g^\sigma)_j \qquad Z_{D(G)}(\alpha_i) = (1 \otimes 1)_i \qquad \sigma \in \{-1, 1\}$$

Using the Wirtinger presentation, the corresponding knot invariant $Z_{D(G)}$ may be interpreted as

$$Z_{D(G)}(K) = \sum_\pi \delta^{\pi(\mu)} \otimes \pi(\lambda)$$

Here the sum is over all representations $\pi$ of the fundamental group of the knot complement into $G$ and $\mu, \lambda$ are the canonical meridian and 0-framed longitude of the long knot.

For example the value of the trefoil knot is

$$Z_{D(G)}(X_{15}^+ X_{62}^+ X_{37}^+ \alpha_4^{-1} /\!/ m^{(1234567)}) = \sum_{a,b,c \in G}(\delta^a \delta^{bcb^{-1}} \delta^{bab(ba)^{-1}} \otimes a^{-3}bac)_1$$

We leave these statements as an exercise to the reader as they are neither difficult nor new [9] and also not the subject of the present paper. In the next section we will work out a more relevant and interesting example in full detail.

Before going into specific examples of snarl algebras it should be mentioned that many such algebras can be obtained from applying the Drinfeld double construction [5]. This includes the finite group example just given. The resulting ribbon Hopf algebras always yield snarl algebras but snarl algebras are a little simpler since we do not require the coalgebra structure. In future work we will comment more on these general constructions.

# 4 Alexander polynomial and Weyl algebra

In this section we introduce the Weyl algebra and show how it is a snarl algebra. The corresponding invariant is the Alexander polynomial and serves as both a special case and a warming up example for the invariant treated in the next section.

The Weyl algebra $\mathcal{W}$ is the algebra generated over the ring $R = \mathbb{Q}(t^{\frac{1}{2}})$ by non-commutative power series in $E, F$ subject to the relation $FE - EF = 1$. Using this relation any element of $\mathcal{W}$ may be written as a sum of alphabetically ordered monomials. Infinite series are always understood in the topology suggested by this alphabetic ordering.

The foundation for our computations is the well-known Weyl commutation relation between exponentials.

**Lemma 1.** *In $\mathcal{W}$ we have the following relation:*

$$\mathbf{e}^{yF}\mathbf{e}^{xE} = \mathbf{e}^{yx}\mathbf{e}^{xE}\mathbf{e}^{yF} \tag{9}$$

*Proof.* This follows directly from the identity $F^b E^a = \sum_j \frac{a!b!E^{a-j}F^{b-j}}{(a-j)!j!(b-j)!}$ that may be proven by induction. $\square$

For convenience we formalize our alphabetic approach to $\mathcal{W}$ a bit more by introducing $\mathcal{V} = R[[e, f]]$. The map $\mathbb{O}: \mathcal{V} \to \mathcal{W}$ given by $\mathbb{O}(e^i f^j) = E^i F^j$ is a bijection, this follows from ordering. We may use it to pull back the product on $\mathcal{W}$ as follows. Define $m(g, g') = \mathbb{O}^{-1}(\mathbb{O}(g)\mathbb{O}(g'))$, for example $m(ef, ef) = e^2 f^2 + ef$. Our convention is that this special product on $\mathcal{V}$ will be denoted $m$, while the ordinary product of power-series is denoted by juxtaposition. Since $\mathcal{W}$ is associative the same is true for $\mathcal{V}$ with the pulled back product $m$.

For later use we also define renaming operations, these are algebra maps $r_j^i : \mathcal{V}^{\otimes S \sqcup \{i\}} \to \mathcal{V}^{\otimes S \sqcup \{j\}}$ by $r_j^i(e_i) = e_j$ and $r_j^i(f_i) = f_j$ and the identity on the other algebra generators. More generally for any ordered partition $\tau$ of $S$ into $b$ parts and $b$-element sequence $B$ we set $r_B^\tau(e_i) = e_{B_j}$ if $i \in \tau^j$ and the same for $f_i$. We sometimes use the short hand $r_k^{ij} = r_k^i /\!/ r_k^j$ and also $r^\tau = r_{(1,2,\ldots,b)}^\tau$.

We now aim to develop techniques for showing $\mathcal{V}$ is a snarl algebra with if we set:

$$Z(X_{ij}^\sigma) = t^{-\frac{\sigma}{2}} \mathbf{e}^{(1-t^\sigma)(e_i - e_j)f_j} \qquad Z(\alpha_i^\sigma) = t^{-\frac{\sigma}{2}}$$

The main difficulty is to find a good formula for the multiplication $m$ on elements such as the above. Let $S$ be a finite set of labels. Lemma 1 tells us that in $\mathcal{V}^{\otimes S}$ we have $m_k^{ij}\mathbf{e}^{bf_i+ae_j} = \mathbf{e}^{ba+ae_k+bf_k}$. Setting $x = (x_s)_{s \in S}$, $y = (y_s)_{s \in S}$ and $e = (e_s)_{s \in S}$ and $f = (f_s)_{s \in S}$ we may write any element of $g \in \mathcal{V}^{\otimes S}$ as $g(e, f) = g(\partial_x, \partial_y)\mathbf{e}^{ex+yf}|_{x=y=0}$. This allows us to compute $m_k^{ij}(g) =$

$$m_k^{ij}(g(\partial_x, \partial_y)\mathbf{e}^{ex+yf}|_{x=y=0}) = g(\partial_x, \partial_y)m_k^{ij}(\mathbf{e}^{ex+yf})|_{x=y=0} = g(\partial_x, \partial_y)\mathbf{e}^{y_i x_j + xe + yf}|_{x=y=0} /\!/ r_k^{ij}$$

More specifically if $g = \mathbf{e}^{eQf}$ for some square matrix $Q$, as is the case for the fundamental snarls then

$$m_k^{ij}(\mathbf{e}^{eQf}) = \mathbf{e}^{\partial_x Q \partial_y}\mathbf{e}^{y_i x_j + xe + yf}|_{x=y=0} /\!/ r_k^{ij}$$

The lemma below will show us how to simplify this.

**Lemma 2.** *Given square matrices $W, Q$ such that $\det(I - WQ) \neq 0$ and vectors $x, y, e, f$ we have*

$$\mathbf{e}^{\partial_x Q \partial_y} \mathbf{e}^{yWx+ex+yf} = \frac{\mathbf{e}^{(eQ+y)(I-WQ)^{-1}(f+Wx)+ex}}{\det(I-WQ)} \qquad (10)$$

*Proof.* We claim that both sides of the equation satisfy the system of differential equations

$$\partial_{Q_{ij}} \Psi = \partial_{x_i} \partial_{y_j} \Psi \qquad \Psi|_{Q=0} = \mathbf{e}^{yWx+ex+yf}$$

There is only one power series in the commuting variables $Q_{i,j}, x_i, y_j$ satisfying these equations. Indeed we may use the differential equation to express the coefficient of any monomial in terms of coefficients of monomials whose joint degree in the $Q$-variables is lower. It thus suffices to prove our claim that both sides satisfy the differential equations. We focus on the right-hand side as the left hand side is clear. Set $A = I - WQ$ so the exponent is $V = (eQ + y)A^{-1}(f + Wx) + ex$. Then

$$\partial_{Q_{ij}} \frac{\mathbf{e}^V}{\det A} = \frac{\mathbf{e}^V}{\det A} \left( \partial_{Q_{ij}} V - \frac{\partial_{Q_{ij}}(\det A)}{\det A} \right)$$

For the first term we use $\partial_{Q_{ij}}(A^{-1})_{rs} = (A^{-1}W)_{ri} A^{-1}_{js}$ to find

$$\partial_{Q_{ij}} V = e_i (A^{-1}(f+Wx))_j + ((eQ+y)A^{-1}W)_i (A^{-1}(f+Wx))_j$$

For the second term we compute $\partial_{Q_{ij}}(\det A) = -(\mathrm{adj}(A)W)_{ji}$ so $-\frac{\partial_{Q_{ij}}(\det A)}{\det A} = (A^{-1}W)_{ji}$ This matches the other side of the equation:

$$\partial_{x_i} \partial_{y_j} \frac{\mathbf{e}^V}{\det A} = \frac{\mathbf{e}^V}{\det A} ((eQ+y)A^{-1}W + e)_i (A^{-1}(f+Wx))_j + (A^{-1}W)_{ji})$$

□

Returning to our computation of $m_k^{ij}(\mathbf{e}^{eQf})$, taking $W = E_{ij}$ gives $(I-WQ)^{-1} = I + \frac{1}{1-Q_{ji}} \sum_b Q_{jb} E_{ib}$ and so $\det(I - WQ) = 1 - Q_{ji}$ and $Q(I-WQ)^{-1} = Q + \frac{1}{1-Q_{ji}} \sum_{ab} Q_{jb} Q_{ai} E_{ab}$. This means that

$$m_k^{ij}(\mathbf{e}^{eQf}) = \frac{\mathbf{e}^{eQf + \frac{1}{1-Q_{ji}} \sum_{ab} e_a Q_{ai} Q_{jb} f_b}}{(1-Q_{ji})} /\!/ r_k^{ij}$$

Sometimes it is convenient to do many multiplications (stitchings) at once. For this we generalize the above discussion to prove

**Lemma 3.** *For any ordered partition of the labels $\tau = (\tau^1, ..., \tau^n)$ and an n-element set of new labels $L$ we may describe $m_L^\tau$ as follows. First define $W = \sum_{\{(i,j) \mid \exists s: i, j \in \tau^s, i \prec j\}} E_j^i$. Here $\prec$ refers to the ordering of the elements in $\tau^k$.*

$$m_L^\tau \left( \frac{\mathbf{e}^{eQf}}{\Delta} \right) = \frac{1}{\Delta \det(I-WQ)} \mathbf{e}^{Q(I-WQ)^{-1}} /\!/ r_L^\tau$$

*For any constant $\Delta$ and matrix $Q$ whose entries are indexed by $L$.*

*Proof.* Without loss of generality we restrict ourselves to the special case $\tau = (1, 2, \ldots n) = S$. For any $g \in \mathcal{V}^{\otimes S}$ we have $m^\tau(g) =$

$$m_1^\tau(g(\partial_x, \partial_y) \mathbf{e}^{ex+yf}|_{x=y=0}) = g(\partial_x, \partial_y) m_1^\tau(\mathbf{e}^{ex+yf})|_{x=y=0} = g(\partial_x, \partial_y) \mathbf{e}^{yWx+xe+yf}|_{x=y=0} /\!/ r_1^\tau$$

Here we used the commutation relation that follows from lemma 1

$$\mathbf{e}^{x_1 E} \mathbf{e}^{y_1 F} \mathbf{e}^{x_2 E} \mathbf{e}^{y_2 F} \ldots \mathbf{e}^{x_n E} \mathbf{e}^{y_n F} = \mathbf{e}^{yWx} \mathbf{e}^{(x_1+\cdots+x_n)E} \mathbf{e}^{(y_1+\cdots+y_n)F}$$

Now if $g = \frac{\mathbf{e}^{eQf}}{\Delta}$ Lemma 2 finishes the proof. □

In case $L = \{1, ..., k\}$ and $\tau$ is a partition of $\{1, ..., n\}$ we may write out the renaming operation $r^\tau$ explicitly. Define $M$ by $M_{ij} = 1$ if $i \in \tau^j$ and zero otherwise. Then

$$m_L^\tau \left( \frac{1}{\Delta} \mathbf{e}^{eQf} \right) = \frac{1}{\Delta \det(I-WQ)} \mathbf{e}^{eM^T Q(I-WQ)^{-1} Mf}$$

## 4.1 The Weyl algebra is a snarl algebra

In this section we prove that $\mathcal{W}$, or rather its normal ordered version $\mathcal{V}$, with the formulas below really satisfies the axioms of a snarl algebra. We will treat the most important cases and write down the main steps in the computations. Readers that would like to see more details are invited to run the computer implementation described in the Appendix.

$$Z(X_{ij}^\sigma) = t^{-\frac{\sigma}{2}} \mathbf{e}^{(1-t^\sigma)(e_i - e_j)f_j} \qquad Z(\alpha_i^\sigma) = t^{-\frac{\sigma}{2}}$$

One of the benefits of extending knot invariants to local objects like snarls is that checking the Reidemeister moves becomes local too. Each of the axioms is a routine calculation using the formulas we derived above (Lemma 3 and the remark coming after it). Here and in the sequel we often write the invariant of a snarl $T$ as $Z(T) = \Delta_T^{-1} \mathbf{e}^{eQ_T f}$ where $\Delta_T, Q_T$ are the constant and matrix defined by this equation.

First we check
$$Z(X_{31}^{\pm} \alpha_2^{\pm 1} /\!/ m^{(123)}) = Z(\alpha_1^0)$$

To evaluate the left hand side we first consider the disjoint union $D = X_{31}^\sigma \sqcup \alpha_2^\sigma$. $Z(D) = \Delta_D^{-1} \mathbf{e}^{eQ_D f}$, where $\Delta_D = t^\sigma$ and $W$ and $Q_D$ are given below:

$$W = \begin{pmatrix} 0 & 1 & 1 \\ 0 & 0 & 1 \\ 0 & 0 & 0 \end{pmatrix} \qquad Q_D = \begin{pmatrix} t^\sigma - 1 & 0 & 0 \\ 0 & 0 & 0 \\ 1 - t^\sigma & 0 & 0 \end{pmatrix} \qquad (I - WQ)^{-1} = \begin{pmatrix} t^{-\sigma} & 0 & 0 \\ t^{-\sigma} - 1 & 1 & 0 \\ 0 & 0 & 1 \end{pmatrix}$$

The last matrix has determinant $t^{-\sigma}$ and $M^T Q(I - WQ)^{-1} M = 0$, where $M^T = (1, 1, 1)$. It follows that $Z(D /\!/ m^{(123)}) = 1 = Z(\alpha_1^0)$. The other two cases of this equation are similar.

Next consider $Z(X_{12}^- X_{34}^+ \alpha_5 \alpha_6^{-1} /\!/ m^{(13)(4526)}) = Z(\alpha_1^0 \alpha_2^0)$. Again we first study the disjoint union $D = X_{12}^- X_{34}^+ \alpha_5 \alpha_6^{-1}$. Then $\Delta_D = 1$ and $Q, W$ are given below:

$$Q_D = \begin{pmatrix} 0 & \frac{t-1}{t} & 0 & 0 & 0 & 0 \\ 0 & -\frac{t-1}{t} & 0 & 0 & 0 & 0 \\ 0 & 0 & 0 & 1-t & 0 & 0 \\ 0 & 0 & 0 & t-1 & 0 & 0 \\ 0 & 0 & 0 & 0 & 0 & 0 \\ 0 & 0 & 0 & 0 & 0 & 0 \end{pmatrix} \qquad W = \begin{pmatrix} 0 & 0 & 1 & 0 & 0 & 0 \\ 0 & 0 & 0 & 0 & 0 & 1 \\ 0 & 0 & 0 & 0 & 0 & 0 \\ 0 & 1 & 0 & 0 & 1 & 1 \\ 0 & 1 & 0 & 0 & 0 & 1 \\ 0 & 0 & 0 & 0 & 0 & 0 \end{pmatrix} \qquad (I - WQ)^{-1} = \begin{pmatrix} 1 & \frac{(t-1)^2}{t} & 0 & 1-t & 0 & 0 \\ 0 & 1 & 0 & 0 & 0 & 0 \\ 0 & 0 & 1 & 0 & 0 & 0 \\ 0 & -\frac{t-1}{t} & 0 & 1 & 0 & 0 \\ 0 & -\frac{t-1}{t} & 0 & 0 & 1 & 0 \\ 0 & 0 & 0 & 0 & 0 & 1 \end{pmatrix}$$

Since $M^T = \begin{pmatrix} 1 & 0 & 1 & 0 & 0 & 0 \\ 0 & 1 & 0 & 1 & 1 & 1 \end{pmatrix}$ we find $M^T Q(I - WQ)^{-1} M = 0$ and $\det(I - WQ)^{-1}$. Therefore $Z(D /\!/ m^{(13)(4526)}) = 1 = Z(\alpha_1^0 \sqcup \alpha_2^0)$. The second equation is proven similarly.

Next consider $Z(X_{12}^1 X_{34}^1 X_{56}^1 /\!/ m^{(13)(25)(46)}) = Z(X_{12}^1 X_{34}^1 X_{56}^1 /\!/ m^{(35)(16)(24)})$. The disjoint union of the left hand side is $D_L = X_{12}^1 X_{34}^1 X_{56}^1$. Then $\Delta_{D_L} = t^{-\frac{3}{2}}$ and $Q_{D_L}, W_L$ are given below:

$$Q_{D_L} = \begin{pmatrix} 0 & 1-t & 0 & 0 & 0 & 0 \\ 0 & t-1 & 0 & 0 & 0 & 0 \\ 0 & 0 & 0 & 1-t & 0 & 0 \\ 0 & 0 & 0 & t-1 & 0 & 0 \\ 0 & 0 & 0 & 0 & 0 & 1-t \\ 0 & 0 & 0 & 0 & 0 & t-1 \end{pmatrix} \qquad W_L = \begin{pmatrix} 0 & 0 & 1 & 0 & 0 & 0 \\ 0 & 0 & 0 & 0 & 1 & 0 \\ 0 & 0 & 0 & 0 & 0 & 0 \\ 0 & 0 & 0 & 0 & 0 & 1 \\ 0 & 0 & 0 & 0 & 0 & 0 \\ 0 & 0 & 0 & 0 & 0 & 0 \end{pmatrix} \qquad (I - W_L Q_{D_L})^{-1} = \begin{pmatrix} 1 & 0 & 0 & 1-t & 0 & -(t-1)^2 \\ 0 & 1 & 0 & 0 & 0 & 1-t \\ 0 & 0 & 1 & 0 & 0 & 0 \\ 0 & 0 & 0 & 1 & 0 & t-1 \\ 0 & 0 & 0 & 0 & 1 & 0 \\ 0 & 0 & 0 & 0 & 0 & 1 \end{pmatrix}$$

So $\det(I - W_L Q_{D_L}) = 1$ and

$$M_L^T = \begin{pmatrix} 1 & 0 & 1 & 0 & 0 & 0 \\ 0 & 1 & 0 & 0 & 1 & 0 \\ 0 & 0 & 0 & 1 & 0 & 1 \end{pmatrix} \qquad M_L^T Q(I - W_L Q_{D_L})^{-1} M_L = \begin{pmatrix} 0 & 1-t & 1-t \\ 0 & t-1 & -(t-1)t \\ 0 & 0 & (t-1)(t+1) \end{pmatrix}$$

Likewise, the disjoint union for the right hand side is $D_R = D_L$ but with different $W_R, M_R$ as shown below

$$W_R = \begin{pmatrix} 0 & 0 & 0 & 0 & 0 & 1 \\ 0 & 0 & 0 & 1 & 0 & 0 \\ 0 & 0 & 0 & 0 & 1 & 0 \\ 0 & 0 & 0 & 0 & 0 & 0 \\ 0 & 0 & 0 & 0 & 0 & 0 \\ 0 & 0 & 0 & 0 & 0 & 0 \end{pmatrix} \qquad M_R^T = \begin{pmatrix} 0 & 0 & 1 & 0 & 1 & 0 \\ 1 & 0 & 0 & 0 & 0 & 1 \\ 0 & 1 & 0 & 1 & 0 & 0 \end{pmatrix} \qquad (I - W_R Q_{D_R})^{-1} = \begin{pmatrix} 0 & 1-t & 1-t & 0 & 0 & 0 \\ 0 & t-1 & -(t-1)t & 0 & 0 & 0 \\ 0 & 0 & (t-1)(t+1) & 0 & 0 & 0 \\ 0 & 0 & 0 & 0 & 0 & 0 \\ 0 & 0 & 0 & 0 & 0 & 0 \\ 0 & 0 & 0 & 0 & 0 & 0 \end{pmatrix}$$

So $\det(I - W_R Q_{D_R}) = 1$ and $M_R^T Q(I - W_R Q_{D_R})^{-1} M_R = \begin{pmatrix} 0 & 1-t & 1-t \\ 0 & t-1 & -(t-1)t \\ 0 & 0 & (t-1)(t+1) \end{pmatrix}$

For the final equation $Z(X_{12}^1) = Z(\alpha_1^1 \alpha_2^1 X_{34}^{\pm} \alpha_5^{-1} \alpha_6^{-1} /\!/ m^{(135)(246)})$ we consider the right-hand side as the stitching of $D = \alpha_1^1 \alpha_2^1 X_{34}^{\pm} \alpha_5^{-1} \alpha_6^{-1}$, with $\Delta_D = t^{-\frac{1}{2}}$ and

$$Q_D = \begin{pmatrix} 0 & 0 & 0 & 0 & 0 & 0 \\ 0 & 0 & 0 & 0 & 0 & 0 \\ 0 & 0 & 0 & 1-t & 0 & 0 \\ 0 & 0 & 0 & t-1 & 0 & 0 \\ 0 & 0 & 0 & 0 & 0 & 0 \\ 0 & 0 & 0 & 0 & 0 & 0 \end{pmatrix} \qquad W = \begin{pmatrix} 0 & 0 & 1 & 0 & 1 & 0 \\ 0 & 0 & 0 & 1 & 0 & 1 \\ 0 & 0 & 0 & 0 & 1 & 0 \\ 0 & 0 & 0 & 0 & 0 & 1 \\ 0 & 0 & 0 & 0 & 0 & 0 \\ 0 & 0 & 0 & 0 & 0 & 0 \end{pmatrix} \qquad (I - WQ_D)^{-1} = \begin{pmatrix} 1 & 0 & 0 & 1-t & 0 & 0 \\ 0 & 1 & 0 & t-1 & 0 & 0 \\ 0 & 0 & 1 & 0 & 0 & 0 \\ 0 & 0 & 0 & 1 & 0 & 0 \\ 0 & 0 & 0 & 0 & 1 & 0 \\ 0 & 0 & 0 & 0 & 0 & 1 \end{pmatrix}$$

This shows that with $M^T = \begin{pmatrix} 1 & 0 & 1 & 0 & 1 & 0 \\ 0 & 1 & 0 & 1 & 0 & 1 \end{pmatrix}$ we have $M^T Q (I - WQ)^{-1} M = \begin{pmatrix} 0 & 1-t \\ 0 & t-1 \end{pmatrix}$ and $\det(I - WQ) = 1$ as required.

This proves that $\mathcal{W}$ is indeed a snarl algebra. Therefore the corresponding invariant $Z_\mathcal{W}$ is independent of the chosen snarl diagram. In the next section we will see that in fact $Z_\mathcal{W}$ is the Alexander polynomial.

## 4.2 Connection to the Burau representation and the Alexander polynomial

For braids we can make the connection with the Burau representation $\beta_t : B_n \to GL(n)$ [3] p.162. Recall on generators $\sigma_k$ we have $\beta_t(\sigma_k) = I - tE_k^k + tE_{k+1}^k + E_k^{k+1} - E_{k+1}^{k+1}$. The special case $\beta_1$ factors through the symmetric group and just gives the permutation matrix induced by the braid. Also recall that the Alexander polynomial may be defined (up to $\pm t^k$) as $\det_1(I - \beta_t(b))$ where the subscript indicates the first row and column of the matrix need to be deleted before taking the determinant.

**Theorem 4.**  1. Suppose $b$ is a braid. Viewed as a snarl we may compute $Z_\mathcal{W}(b) = t^{-\frac{w}{2}} \mathbf{e}^{e(\beta_t(b)^T \beta_1(b) - I)f}$, where $w$ is the signed sum of the crossings (writhe).

  2. For any knot $K$, viewed as a snarl, $Z_\mathcal{W}(K) = \Delta_K^{-1}$ where $\Delta_K$ is the Alexander polynomial of the knot.

*Proof.* Part 1) Building the braid from a disjoint union $D$ of crossings we find that $\Delta_D = t^{\frac{w}{2}}$. If we label the pieces of $D$ so that labels ending up in the same component are enumerated in order of appearance then $(1 - WQ)$ is upper triangular with ones on the diagonal. This proves $\Delta_b = t^{\frac{w}{2}}$. Moving on to $Q_b$ we aim to prove that
$$\beta_t(b) = \beta_1(b)(I + Q_b)^T \tag{11}$$
Notice that this formula is correct when $b = \sigma_k^{\pm 1}$ is a generator of the braid group or the identity. To prove the general case it suffices to show that the right hand side is multiplicative in $b$. In other words we should show that $\beta_1(b')(I + Q_{b'})^T \beta_1(b)(I + Q_b)^T = \beta_1(b'b)(I + Q_{b'b})^T$. This is equivalent to showing
$$Q_{b'b} = Q_b + Q_b \beta_1(b)^T Q_{b'} \beta_1(b) + \beta_1(b)^T Q_{b'} \beta_1(b) \tag{12}$$

This formula follows directly from stitching the ends of the disjoint union $D = b \sqcup b'$ where we number the components of $b$ by $1, 2, ..n$ as they appear at the bottom and likewise $b'$ by $n + 1, ..2n$. The precise stitching rule is determined by the permutation $\beta_1(b)$ induced by $b$: the end $i$ of $b$ is stitched to the beginning of component $n + \beta_1(i)$ of $b'$. Therefore in the stitching formula the matrices $Q_D, W, M$ have the following block shapes:
$$W = \begin{pmatrix} 0 & \beta_1(b)^T \\ 0 & 0 \end{pmatrix} \quad Q_D = \begin{pmatrix} Q_b & 0 \\ 0 & Q_{b'} \end{pmatrix} \quad M = \begin{pmatrix} I \\ \beta_1(b) \end{pmatrix}$$
The result (12) now follows from computing $Q_{b'b} = M^T Q (1 - WQ)^{-1} M$. Note the inverse is easy to compute as it is an upper-triangular matrix.

Part 2) To see that $Q_K = 0$ for any knot, consider stitching the knot from a disjoint union $D$ of crossings and $\alpha$'s and notice that $M^T Q_D = 0$, hence $Q_K = 0$. Now view the long knot $K$ as the partial closure of braid $b$ to interpret the inverse of the constant as the Alexander polynomial. By conjugating our braid we may assume that only the first top and bottom strands are open and $\beta_1(b)$ is the permutation matrix of the cycle $(1, 2, 3, ..., n)$. At first let us ignore the $\alpha$'s and stitch the ends of the braid directly.

We aim to show that the determinant $\det(I - WQ_b)$, arising from stitching the braid as above, equals the determinant $\det_1(I - \beta_t(b))$ defining the Alexander polynomial. As a first step notice that the bottom row of $WQ_b$ consists of zeroes. Hence $\det(I - WQ_b) = \det_n(I - WQ_b)$, where $\det_k$ is the determinant after deleting the $k$-th row and column.

The formula for $Q_b$ from part 1) indicates proving $\det_1(I - \beta_t(b)) = \det_n(I - W(\beta_t(b)^T \beta_1(b) - I))$ is enough. To prove this identity we note that the right hand side equals

$$\det_n((I - \beta_1(b)^T)(I - W(\beta_t(b)^T \beta_1(b) - I))) = \det_n(I - \beta_1(b)^T - \beta_1(b)^T(\beta_t(b)^T \beta_1(b) - I)) = \det_n(I - \beta_1(b)^T \beta_t(b)^T \beta_1(b))$$

Conjugation by $\beta_1(b)$ and transposition turns it into the desired $\det_1(I - \beta_t(b))$.

Had we included the $\alpha$'s into the stitching the only thing that changes is we pick up a power $t^{-\frac{n}{2}}$ where $n$ is the number of strands in our braid (we are doing a closure to the right so we only pick up $n$ copies of $\alpha^{-1}$). □

In fact the normalization of the Alexander polynomial we compute appears to be symmetric with respect to $t \mapsto t^{-1}$.

# 5 Generalization to the $q$-Weyl algebra

In this section we present the knot invariant built on a deformation of the Weyl algebra. The theory is developed in complete analogy to the case of the usual Weyl algebra and the Alexander polynomial described above. We will also indicate how this invariant relates to the knot invariant $Z_0$ from the introduction.

Define the $q$-Weyl[5] algebra $\mathcal{W}_q$ to be the associative algebra with generators $1, E, F$ and relation $FE = qEF + 1$. To stay as close to the Alexander polynomial as possible we restrict ourselves to the special case where our variable $q$ satisfies $q = 1 + \epsilon$ and $\epsilon^2 = 0$. This means we consider our algebra over the base ring $R = \mathbb{Q}(t^{\frac{1}{2}})[\epsilon]/(\epsilon^2)$. Many of the techniques below apply to more general values of $q$ but we leave this for future work.

Alphabetic ordering of the monomials is used precisely as in the case of $\mathcal{W}$ to deal with infinite series. The commutation relation between exponentials is as follows (recall that $\mathbf{e}^{\epsilon x} = 1 + \epsilon x$).

**Lemma 4.**
$$\mathbf{e}^{yF}\mathbf{e}^{xE} = \mathbf{e}^{xE}\mathbf{e}^{xy}\mathbf{e}^{\epsilon(\frac{y}{2}+E)(\frac{x}{2}+F)}\mathbf{e}^{yF}$$

And more generally:
$$\prod_{i=1}^{n}\mathbf{e}^{x_iE}\mathbf{e}^{y_iF} = \mathbf{e}^{(x_1+..+x_n)E}\mathbf{e}^{\sum_j x_j(y_1+..+y_{j-1})}\mathbf{e}^{\epsilon(\frac{y_1+..y_{j-1}}{2}+E)(\frac{x_j}{2}+x_{j+1}+..x_n+F)}\mathbf{e}^{(y_1+..+y_n)F}$$

*Proof.* If we define $[k] = \frac{1-q^k}{1-q} = k(1+\epsilon\frac{k-1}{2})$ and $[a]! = a!(1+\epsilon\frac{a(a-1)}{4})$ induction as in Lemma 1 shows
$$F^bE^a = \sum_j \frac{[a]![b]!E^{a-j}(1+j\epsilon EF)F^{b-j}}{[a-j]![j]![b-j]!} = \sum_j \frac{a!b!E^{a-j}(1+\epsilon j(\frac{a+b}{2}-1-\frac{3}{4}(j-1))+EF)F^{b-j}}{(a-j)!j!(b-j)!}$$

The relations between the exponentials are a direct consequence. □

We prefer to work in commutative setting so define $\mathcal{V}_q$ to be the power series ring $R[[e,f]]$ as in the Alexander case. The map $\mathbb{O} : \mathcal{V}_q \to \mathcal{W}_q$ defined by $\mathbb{O}(e^a f^b) = E^a F^b$ is used to pull back the multiplication $m$ from $\mathcal{W}_q$ to $\mathcal{V}_q$. Notice $\mathbb{O}$ remains a bijection in this more general set up. We claim the following makes $\mathcal{V}_q$ into a snarl algebra.

$$Z(X_{i,j}^\sigma) = t^{-\frac{\sigma}{2}}\mathbf{e}^{(1-t^\sigma)(e_i-e_j)f_j+\epsilon P} \quad Z(\alpha_i^\sigma) = t^{-\frac{\sigma}{2}}\mathbf{e}^{\epsilon\sigma e_i f_i}$$

where
$$P = \sigma\left(\left(\frac{1-t^\sigma}{4}e_if_j\right)^2 - \left(\frac{1+t^\sigma}{4}e_jf_j\right)^2 + t^\sigma e_if_j(e_jf_i + \frac{1-t}{2}((\sigma+1)e_jf_j + (\sigma-1)e_if_i))\right)$$

To work with such formulas we need to develop good stitching (multiplication) formulas as we did in the Alexander case. Since the arguments are entirely analogous we summarize the results in the following lemma.

**Lemma 5.** *For any ordered partition of the labels $\tau = (\tau^1,..,\tau^n)$ and an $n$-element set of new labels $L$, define $W = \sum_{\{(i,j)|\exists s: i,j\in\tau^s, i\prec j\}} E_j^i$ and $A = (I-WQ)^{-1}$. Here $\prec$ refers to the ordering of the elements in $\tau^k$. Also set*

$$S(x,y) = \frac{t+1}{t-1}\sum_{m=1}^n \sum_{j\in\tau^m}\left(\sum_{i\prec j\in\tau^m} y_i\right) x_j \left(\frac{1}{2}\sum_{k\prec j\in\tau^m} y_k + \tilde{e}_j\right)\left(\frac{x_j}{2} + \sum_{j\prec l\in\tau^m} x_l + \tilde{f}_j\right)$$

Then
$$m_L^\tau\left(\frac{\mathbf{e}^{eQf+\epsilon P}}{\Delta}\right) = (1+\epsilon P(\partial_x,\partial_y) + \epsilon S(\partial_e,\partial_f))\frac{\mathbf{e}^{(eQ+y)A(f+Wx)+ex}}{\Delta\det(I-WQ)}\Big|_{x=y=0;\tilde{e},\tilde{f}\mapsto e,f}/\!/r_L^\tau$$

*For any constant $\Delta$ and matrix $Q$ whose entries are indexed by $L$.*

*Proof.* The proof is analogous to the $\mathcal{W}$ case. Without loss of generality we consider the case $\tau = (1,2\ldots,n)$. For any $g \in \mathcal{V}_q^{\otimes S}$ we have $m^\tau(g) = $

$$m^\tau(g(\partial_x,\partial_y)\mathbf{e}^{ex+yf}|_{x=y=0}) = g(\partial_x,\partial_y)m^\tau(\mathbf{e}^{ex+yf})|_{x=y=0} = g(\partial_x,\partial_y)\mathbf{e}^{yWx+xe+yf+\epsilon S(x,y)}|_{x=y=0}/\!/r^\tau$$

Here we used the commutation relation from lemma 4. Now set $g = \Delta^{-1}\mathbf{e}^{eQf+\epsilon P(e,f)}$ and recall $\epsilon^2 = 0$ so that we find

$$m^\tau(g) = \Delta^{-1}\mathbf{e}^{\partial_x Q\partial_y + \epsilon P(\partial_x,\partial_y)}\mathbf{e}^{yWx+xe+yf+\epsilon S(x,y)}|_{x=y=0;\tilde{e},\tilde{f}\mapsto e,f}/\!/r^\tau =$$
$$(1+\epsilon P(\partial_x,\partial_y) + \epsilon S(\partial_e,\partial_f))\Delta^{-1}\mathbf{e}^{\partial_x Q\partial_y}\mathbf{e}^{yWx+xe+yf}|_{x=y=0;\tilde{e},\tilde{f}\mapsto e,f}/\!/r^\tau$$

Applying Lemma 2 gives us

$$(1+\epsilon P(\partial_x,\partial_y))(1+\epsilon S(\partial_e,\partial_f))\frac{\mathbf{e}^{(eQ+y)(I-WQ)^{-1}(f+Wx)+ex}}{\Delta\det(I-WQ)}\Big|_{x=y=0;\tilde{e},\tilde{f}\mapsto e,f}/\!/r^\tau$$

as required. □

---
[5]Often called q-Heisenberg algebra [8].

The above formula may be simplified slightly since both $P$ and $S$ only depend on two of the four sets of variables. We may rewrite it as $m_L^\tau(\frac{\mathbf{e}^{eQf+\epsilon P}}{\Delta}) = \frac{Z_C + \epsilon Z_S + \epsilon Z_P}{\Delta \det(I-WQ)} /\!/ r_L^\tau$ where

$$Z_C = \mathbf{e}^{eQAf} \quad Z_S = S(\partial_e, \partial_f)\mathbf{e}^{eQAf}|_{\tilde{e},\tilde{f}\mapsto e,f} \quad Z_P = P(\partial_x, \partial_y)\mathbf{e}^{e(QAW+I)x+yA(f+Wx)}|_{x=y=0}$$

Now that we know how to multiply (stitch) exponentials in $\mathcal{W}_q$ (or rather $\mathcal{V}_q$) we are in a position to turn it into a snarl algebra.

**Proposition 1.** $\mathcal{W}_q$ is a snarl algebra when we set

$$Z_{\mathcal{W}_q}(X_{i,j}^\sigma) = t^{-\frac{\sigma}{2}}\mathbf{e}^{(1-t^\sigma)(e_i-e_j)f_j+\epsilon P} \quad Z_{\mathcal{W}_q}(\alpha_i^\sigma) = t^{-\frac{\sigma}{2}}\mathbf{e}^{\epsilon\sigma e_i f_i}$$

where

$$P = \sigma\left(\left(\frac{1-t^\sigma}{4}e_i f_j\right)^2 - \left(\frac{1+t^\sigma}{4}e_j f_j\right)^2 + t^\sigma e_i f_j(e_j f_i + \frac{1-t}{2}((\sigma+1)e_j f_j + (\sigma-1)e_i f_i))\right)$$

*Proof.* In checking the snarl algebra axioms we only need to pay attention to the $\epsilon$-dependent part as the rest together with the matrices $Q, W$ etc are already computed in Section 4.1 where we proved that $\mathcal{W}$ was a snarl algebra with a compatible snarl structure. Given the formulas for stitching in Lemma 5 checking the axioms is a straightforward if tedious calculation. As such routine calculations are best done by computer we refer the reader to the Mathematica implementation in the Appendix for more details. □

By section any snarl algebra gives rise to a knot invariant. Hence we now have a knot invariant $Z_{\mathcal{W}_q}$ coming from the snarl algebra $\mathcal{W}_q$. In the next section we explore some of its properties. We end this section by showing how the formula for $Z_0$ from the introduction follows from the above definition of $Z_{\mathcal{W}_q}$.

**Proposition 2.** For any knot $K$, let $C_K$ be the coefficient of $\epsilon$ of the constant part of $Z_{\mathcal{W}_q}$, i.e. when we set $e = f = 0$. We have $C_K = Z_0/\Delta_K^2$, where $Z_0$ was defined in the introduction.

*Proof.* Our set up is that of a knot diagram with $n$ crossings and cuaps numbered in order of appearance as described in the introduction. Denote the disjoint union of these elements as $D$. The knot $K$ is obtained as $K = D /\!/ m_1^{(12\ldots n)}$. To match the description of $Z_0$ we will follow the instructions for computing the $\epsilon$-part of $Z_{\mathcal{W}_q}$ cutting a few corners along the way.

We start with the matrix $W = \sum_{i<j} E_j^i$. Next, the matrix $Q_D$ is $Q_D = \sum_{X_{i,j}^s}(1-t^s)(E_j^i - E_j^j)$. Set $A = (1-WQ)^{-1}$ and $S$ is as above. Since we are only interested in the constant term $C_K$ we may ignore all parts that only contribute $e$ or $f$. The formula becomes:

$$C_K = P_D(\partial_x, \partial_y)\mathbf{e}^{yAWx}|_{x=y=0} + \bar{S}(\partial_e, \partial_f)\mathbf{e}^{eQAf}|_{e=f=0}$$

Here $\bar{S}$ is a simplified version of $S$ without any irrelevant terms such as $\tilde{e}_j$ and $\tilde{f}_j$:

$$\bar{S} = \frac{t+1}{t-1}\sum_{m=1}^n \sum_{j=2}^n \frac{1}{2}\left(\sum_{i,k<j} y_i y_k\right)\left(\frac{x_j^2}{2} + \sum_{j<l} x_j x_l\right)$$

The first term in $C_K$ is a sum over all crossings and cuaps where each monomial $e_i f_j e_k f_l$ contributes $(AW)_j^i (AW)_l^k + (AW)_l^i (AW)_j^k$ and $x_i y_j$ contributes $(AW)_j^i$. Likewise each monomial $y_i x_j y_k x_l$ of $\bar{S}$ contributes $(QA)_j^i (QA)_l^k + (QA)_l^i (QA)_j^k$. As these expressions are nearly homogeneous in $A$ and $Q$ we chose to rescale them: $Q$ is replaced by $\bar{Q} = Q/(t^{\frac{1}{2}} - t^{-\frac{1}{2}})$ and we set $B = A^{-1} = I - (t^{\frac{1}{2}} - t^{-\frac{1}{2}})W\bar{Q}$ as in the introduction. In terms of $B$ we have $A = \mathrm{adj}(B)/\det(B)$. The reader should now recognize the expressions for $Z_H$ and $Z_G$ as the contributions of the crossings/cuaps and the polynomial $\bar{S}$ scaled suitably with $H = \mathrm{adj}(B)W$ and $G = Q\mathrm{adj}(B)$ as described.

We already know that $\Delta_K^{-1} = c^{\frac{1}{2}}\det(B)$ is the normalized Alexander polynomial, where $c$ is the correction that comes from $\Delta_D$. This prompts us to scale $P_K$ by $\Delta_K^2$ and not just by $\det(B)^2$ as we were about to. □

## 6 Properties of $Z_{\mathcal{W}_q}$

In this section we list some (conjectured) properties of the invariant $Z_{\mathcal{W}_q}$ coming from the $q$-Weyl algebra. For simplicity we restrict ourselves to the case for knots.

**Theorem 5.** *If a knot diagram for $K$ is stitched from $n$ fundamental snarls ($X$ and $\alpha$) then $\mathcal{O}(n^6)$ operations in the ring $\mathbb{Z}[t^{\frac{1}{2}}, t^{-\frac{1}{2}}]$ suffice to compute $Z_{\mathcal{W}_q}$.*

*Proof.* By proposition 2 it suffices to consider computing $Z_0$ as described in the introduction. Computing the matrices $G$ and $H$ may be done in $\mathcal{O}(n^4)$ steps. Their matrix entries are Laurent polynomials of degree $\mathcal{O}(n)$ with coefficients bounded by $\mathcal{O}(2^n)$. Since $Z_G$ has $\mathcal{O}(n^4)$ terms and each term consists of multiplying two entries of $G$, computing $Z_G$ takes $\mathcal{O}(n^6)$ operations in $\mathbb{Z}[t^{\frac{1}{2}}, t^{-\frac{1}{2}}]$. Computing $Z_H$ is similar but faster as it involves less terms. □

The above bound can be improved on quite a bit. For example a divide and conquer approach as in [2] is expected to bring the complexity down to $\mathcal{O}(n^5)$. For general snarl diagrams Lemma 5 and similar arguments as above predict the number of operations is $\mathcal{O}(n^8)$.

We remark that the algorithm being polynomial time is a qualitative feature that is independent of the precise notion of complexity of the input used. This is because converting our knot to any other reasonable format can be done in polynomial time (usually linear time). The precise cost of the ring-operations needed is also of little consequence to polynomiality as all such operations can be done in polynomial time.

We end this section with three conjectures on our invariant $Z_{\mathcal{W}_q}(K) = \frac{1}{\Delta_K}\mathbf{e}^{\epsilon Q_K f + \epsilon P_K}$ or rather the normalization

$$\rho_1(K;t) = -\frac{t\Delta_K^2}{(1-t)^2}\left(P_K - t\frac{\mathrm{d}}{\mathrm{d}t}\log \Delta_K\right)$$

**Conjecture 1.** *The normalized knot invariant $\rho_1(K;t)$ is a Laurent polynomial in $t$ with the following symmetries.*

1. *If $-K$ denotes the mirror image of a knot $K$ then we have $\rho_1(-K) = -\rho_1(K)$.*
2. *$\rho_1(K)$ is symmetric with respect to $t \mapsto t^{-1}$.*
3. *$\rho_1(K)$ is invariant under reversing the orientation of $K$.*

In particular, $\rho_1$ vanishes on amphicheiral knots as conjectured by Rozansky for his related invariant, see below. The conjecture was checked experimentally for all knots up to 12 crossings.

As an illustration we list the value of $\rho_1$ on the family of alternating torus knots $T(2, \pm(2p+1))$ where $p \in \mathbb{Z}_{>0}$ is the closure of the braid $\sigma_1^{\mp(2p+1)}$ in $B_2$.

$$\rho_1 T(2,\pm(2p+1)) = \mp\sum_{k=0}^{p-1}\frac{1}{2}(p-k)(p+k+1)(t^{2k+1} - t^{-2k-1})$$

Recall that the Alexander polynomial bounds the genus $g$ as follows. Denote by maxdeg $f(t)$ the highest power of $t$ in Laurent polynomial $f$. When written in its symmetric form, the Alexander polynomial satisfies maxdeg $\Delta \leq g$. A similar but sometimes more powerful bound is provided by $\rho_1$.

**Conjecture 2.** *For a knot $K$ with genus $g$ we have $\mathrm{maxdeg}(\rho_1) \leq 2g - 1$.*

This was checked experimentally for all knots up to 12-crossings and sometimes improves the genus bound given by the Alexander polynomial. For example the knot $12_{313}$ has genus 2 [4] and trivial Alexander polynomial, but the maximal power in $t$ of $\rho_1$ is 2. We expect this conjecture to follow from a formula for $Z$ in terms of a Seifert surface for the knot.

**Conjecture 3.** *$\rho_1(K;t^2) = \frac{t^2}{(1-t^2)^2}P^{(1)}(K;t)$ where $P^{(1)}$ is the invariant of Rozansky defined in [14].*

This conjecture is true for all knots up to nine crossings for which $P^{(1)}$ was computed in [12]. Briefly, Rozansky showed that the colored Jones polynomial may be expanded in $h = q - 1$ as

$$J_\alpha(q) = \sum_{n\geq 0} h^n(\sum_{0\leq m\leq n} D_{m,n}(\alpha h)^{2m})$$

such that the coefficients have the property $\sum_{m\geq 0} D_{m,n+2m}(\alpha h)^{2m} = \frac{P^{(n)}(K;q^\alpha)}{\Delta_K^{2n+1}(q^{\frac{\alpha}{2}} - q^{-\frac{\alpha}{2}})}$ where $\Delta_K$ is as above and $P^{(n)}(K;t)$ is a Laurent polynomial.

The invariant $P^{(1)}$ is also known as the 2-loop invariant and was studied from the Kontsevich integral point of view by Ohtsuki in [10]. His work includes a similar genus bound so that the second conjecture will follow from the last.

Instead of deriving our invariant as an expansion of the colored Jones polynomial, we in some sense expanded the underlying quantum group itself. Our invariant may be interpreted as the universal invariant [11] for some simplified versions of $U_q(sl_2)$ or rather the Drinfeld double [5] of the universal enveloping algebra of the Borel subgroup of $sl_2$. In a sequel to this paper we will explain a general theory for solvable approximation of Lie algebras and how it may be used to derive effective knot invariants. In particular we expect polynomial time algorithms for all of Rozansky's invariants $P^{(n)}$. Our approach is also expected to yield formulas for strand reversal and doubling using Hopf-algebraic techniques. This might yield a proof of the first two conjectures.

## 7 Summary and outlook

After formulating a convenient local notion of knot diagrams called snarl diagrams, we constructed a powerful new knot invariant from the $q$-Weyl algebra. This algebra has two generators and one relation $EF - qFE = 1$. We considered the cases $q = 1$ and $q = 1 + \epsilon$ with $\epsilon^2 = 1$. The first yields the Alexander polynomial, while the latter is new. Both may be computed in polynomial time using normally ordered exponentials. As such our invariant is the strongest known knot invariant computable in polynomial time, this is illustrated in the table in the appendix: all knots up to ten crossings are distinguished.

In future work we plan to show that our techniques apply to algebras much more general than the Weyl algebras presented here. For example we expect that the invariants coming from appropriately truncated Drinfeld doubles of universal enveloping algebras of solvable Lie algebras yields similar invariants. All computable in polynomial time. In this context the $q$-Weyl algebra arose from considering the two-dimensional non-abelian Lie algebra. We also hope to clarify connections to classical invariants such as the genus and with Rozansky's work on the colored Jones polynomial.

DEPARTMENT OF MATHEMATICS, UNIVERSITY OF TORONTO, TORONTO ONTARIO M5S 2E4, CANADA
*E-mail address:* drorbn@math.toronto.edu
*URL:* http://www.math.toronto.edu/drorbn

MATHEMATISCH INSTITUUT, UNIVERSITEIT LEIDEN, NIELS BOHRWEG 1, 2333 CA LEIDEN, HOLLAND
*E-mail address:* r.i.van.der.veen@math.leidenuniv.nl
*URL:* http://www.rolandvdv.nl

# Appendix: Implementation

In this appendix we briefly describe an implementation of the knot invariant described in this paper in Mathematica. This serves two purposes. First it allows us to automatically verify the snarl axioms hold as claimed in Proposition 1. It also allows us to compile a table of the values of our invariant on the table of prime knots up to 10 crossings. The program and table are available at http://www.rolandvdv.nl/MLA/

The program encodes our invariant $Z_{\mathcal{W}_q}(K) = \frac{1}{\Delta} e^{eQf + \epsilon P}$ as $\mathbb{E}[\Delta, Q, P]$, where $\Delta$ is a Laurent polynomial in $t$ and $Q$ (really $eQf$) is a quadratic in $e_i, f_i$ with coefficients in $\mathbb{Q}(t^{\frac{1}{2}})$ and $P$ is a quartic in the same variables.

To specify the stitching $m^\tau$ we use a list of disjoint subsets $\tau^i$ of the label set. For simplicity, after stitching all components in $\tau^i$ are renamed $i$. To improve clarity, disjoint union is written as juxtaposition.

### The program

Introducing the canonical form, disjoint union and some utilities:

```
CF[E[Δ_, Q_, P_]] := Simplify /@ E[Δ, Q, P];
E /: E[Δ1_, Q1_, P1_] E[Δ2_, Q2_, P2_] := E[Δ1 Δ2, Q1 + Q2, P1 + P2];
Pol2Mat[P_, L_] := Table[Coefficient[P, e_{L[[i]]} f_{L[[j]]}], {i, Length[L]}, {j, Length[L]}]
Mat2Pol[M_, L_] := Table[e_{L[[i]]}, {i, Length[L]}].M.Table[f_{L[[i]]}, {i, Length[L]}]
```

The program for stitching, implementing Lemma 5:

```
m_τ_[E[Δ_, Q_, P_]] := Block[{L, n, T, Po, W, A, S, exp, newP},
  {L = Flatten[τ], n = Length[L];, T[s_] := Table[s_i, {i, L}], Po_{i_,j_} := Position[τ, i][[1, j]]};
  W = Sum[If[Po_{i,1} == Po_{j,1} && Po_{i,2} < Po_{j,2}, e_i f_j, 0], {i, L}, {j, L}];
  A = Mat2Pol[Inverse[IdentityMatrix[n] - Pol2Mat[W, L].Pol2Mat[Q, L]], L];
  Fi[j_, k_] := Take[τ[[k]], Position[τ[[k]], j][[1, 1]] - 1];
  La[j_, k_] := Take[τ[[k]], {Position[τ[[k]], j][[1, 1]] + 1, Length[τ[[k]]]}];
  S = (t+1)/(t-1) Sum_{k=1}^{Length[τ]} Sum[Sum[y_i, {i, Fi[j, k]}] x_j (Sum[1/2 y_i, {i, Fi[j, k]}] + ẽ_j)
        (1/2 x_j + Sum[x_i, {i, La[j, k]}] + f̃_j), {j, τ[[k]]}];
  exp = e^((T[e].Pol2Mat[Q,L]+T[y]).Pol2Mat[A,L].(T[f]+Pol2Mat[W,L].T[x])+T[e].T[x]);
  E[Δ / Det[Pol2Mat[A, L]], Mat2Pol[Pol2Mat[Q, L].Pol2Mat[A, L], L],
    exp^-1 Total[CoefficientRules[P + S, Join[T[e], T[f], T[x], T[y]]] /.
      (L_ → c_) ⧴ (c D @@ Prepend[Partition[Riffle[Join[T[x], T[y], T[e], T[f]], L], 2], exp])]
    /. {ẽ_i_ → e_i, f̃_i_ → f_i, x_ → 0, y_ → 0}] /. {z_{-i_} ⧴ z_{Po_{i,1}}} // CF]
```

The values of $\mathbb{Z}_{W_q}$ on the fundamental tangles:

```
R^s_{i_,j_} := E[t^{-s/2}, (1 - t^s) (e_i - e_j) f_j, s (((1-t^s)/2 e_i f_j)^2 - ((1+t^s)/2 e_j f_j)^2 +
    t^s e_i f_j (e_j f_i + (1-t)/2 ((s+1) e_j f_j + (s-1) e_i f_i)))] // CF

u^s_{i_} := E[t^{-s/2}, 0, s e_i f_i]
```

## Verification of Proposition 1

The following commands verify the computations required for Proposition 1 and also Section 4.1, the output is trivial as expected.

$$\{\overset{1}{R}_{1,2}\overset{-1}{u}_3 \mathbin{/\!/} m_{\{\{1,3,2\}\}},\ \overset{1}{R}_{1,2} u_{1,3} \mathbin{/\!/} m_{\{\{2,3,1\}\}},\ \overset{-1}{R}_{1,2}\overset{1}{u}_3 \mathbin{/\!/} m_{\{\{1,3,2\}\}},\ \overset{-1}{R}_{1,2}\overset{-1}{u}_3 \mathbin{/\!/} m_{\{\{2,3,1\}\}}\}$$

$$\{\overset{-1}{R}_{1,2}\overset{1}{R}_{3,4}\overset{1}{u}_5\overset{-1}{u}_6 \mathbin{/\!/} m_{\{\{1,3\},\{4,5,2,6\}\}},\ \overset{1}{R}_{1,2}\overset{-1}{R}_{3,4}\overset{1}{u}_5\overset{-1}{u}_6 \mathbin{/\!/} m_{\{\{5,1,6,3\},\{4,2\}\}}\}$$

$$\{\left(\overset{1}{R}_{1,2}\overset{1}{R}_{3,4}\overset{1}{R}_{5,6} \mathbin{/\!/} m_{\{\{1,3\},\{2,5\},\{4,6\}\}}\right) - \left(\overset{1}{R}_{1,2}\overset{1}{R}_{3,4}\overset{1}{R}_{5,6} \mathbin{/\!/} m_{\{\{3,5\},\{1,6\},\{2,4\}\}}\right),$$

$$\left(\overset{-1}{R}_{1,2}\overset{-1}{R}_{3,4}\overset{-1}{R}_{5,6} \mathbin{/\!/} m_{\{\{1,3\},\{2,5\},\{4,6\}\}}\right) - \left(\overset{-1}{R}_{1,2}\overset{-1}{R}_{3,4}\overset{-1}{R}_{5,6} \mathbin{/\!/} m_{\{\{3,5\},\{1,6\},\{2,4\}\}}\right)\}$$

$$\{\left(\overset{1}{u}_1\overset{1}{u}_2\overset{1}{R}_{3,4}\overset{-1}{u}_5\overset{-1}{u}_6 \mathbin{/\!/} m_{\{\{1,3,5\},\{2,4,6\}\}}\right) - \overset{1}{R}_{1,2},\ \left(\overset{1}{u}_1\overset{1}{u}_2\overset{-1}{R}_{3,4}\overset{-1}{u}_5\overset{-1}{u}_6 \mathbin{/\!/} m_{\{\{1,3,5\},\{2,4,6\}\}}\right) - \overset{-1}{R}_{1,2},$$

$$\left(\overset{-1}{u}_1\overset{-1}{u}_2\overset{1}{R}_{3,4}\overset{1}{u}_5\overset{1}{u}_6 \mathbin{/\!/} m_{\{\{1,3,5\},\{2,4,6\}\}}\right) - \overset{1}{R}_{1,2},\ \left(\overset{-1}{u}_1\overset{-1}{u}_2\overset{-1}{R}_{3,4}\overset{1}{u}_5\overset{1}{u}_6 \mathbin{/\!/} m_{\{\{1,3,5\},\{2,4,6\}\}}\right) - \overset{-1}{R}_{1,2}\}$$

{E[1, 0, 0], E[1, 0, 0], E[1, 0, 0], E[1, 0, 0]}

{E[1, 0, 0], E[1, 0, 0]}

{0, 0}

{0, 0, 0, 0}

## Sample output

As a sample calculation we apply our algorithm to a snarl diagram of the figure eight knot. The disjoint union of the fundamental snarls it is built up from is shown in figure 4.

$$\overset{1}{R}_{1,2}\overset{1}{R}_{3,4}\overset{-1}{R}_{5,6}\overset{-1}{R}_{7,8}\overset{-1}{u}_9\overset{1}{u}_{10} \mathbin{/\!/} m_{\{\{2,3,8,9,5,4,10,1,6,7\}\}}$$

$$E\left[3 - \frac{1}{t} - t,\ 0,\ \frac{-1+t^2}{1-3t+t^2}\right]$$

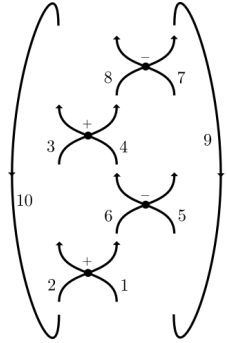

Figure 4: A diagram for the figure eight knot $4_1$ about to be stitched together from a disjoint union of fundamental snarls. Rotation numbers $\rho$ are not listed but should be clear from the picture.